\documentclass[12pt]{amsart}
\usepackage{amsmath}
\usepackage{amssymb}
\usepackage{latexsym}

\newcommand{\Zint}{{\bf {Z}}}    
     
\newcommand{\Rea}{{\bf {R}}}      
\newcommand{\Cplx}{{\bf {C}}}     

\newtheorem{prop}{\bf Proposition}
\newtheorem{thm}[prop]{\bf Theorem}

\newtheorem{cor}[prop]{\bf Corollary}

\begin{document}

\title[Integral transformation and Darboux transformation]{Integral transformation and Darboux transformation of Heun's differential equation}


\author{Kouichi Takemura}
\address{Department of Mathematical Sciences, Yokohama City University, 22-2 Seto, Kanazawa-ku, Yokohama 236-0027, Japan.}
\email{takemura@yokohama-cu.ac.jp}

%

\begin{abstract}
We review Darboux-Crum transformation of Heun's differential equation.
By rewriting an integral transformation of Heun's differential equation into a form of elliptic functions, we see that the integral representation is a generalization of Darboux-Crum transformation.
We also consider conservation of monodromy with respect to the transformations.
\end{abstract}

\maketitle


\section{Introduction}

Heun's differential equation is defined by
\begin{equation}
\left\{ \left(\frac{d}{dz}\right) ^2 \! + \left( \frac{\gamma}{z}+\frac{\delta }{z-1}+\frac{\epsilon}{z-t}\right) \frac{d}{dz} \right.\\
\left. +\frac{\alpha \beta z -q}{z(z-1)(z-t)} \right\} y=0,
\label{eq:Heun}
\end{equation}
with the condition $\gamma +\delta +\epsilon =\alpha +\beta +1$ (\cite{Ron,SL}).
It has four singularities $ \{ 0, 1, t, \infty \}$ on the Riemann sphere $\Cplx \cup \{ \infty \}$ and all the singularities are regular.
Heun's differential equation is a standard form of a second-order Fuchsian differential equation with four singularities.
Heun's differential equation and the differential equations of its confluent type appear in several systems of physics including analysis of black holes, AdS/CFT correspondence (\cite{KSY}), crystal transition (\cite{SL}), fluid dynamics (\cite{CH}).
A standard form of a second-order Fuchsian differential equation with three singularities is the hypergeometric differential equation, which is celebrated in physics and mathematics.
Although analysis of Heun's differential equation is much difficult than that of the hypergeometric equation, several solutions are known.
For example, Heun polynomials (\cite{Ron}) are polynomial solutions of Heun's differential equation, which are related with quasi-exact solvability.
Finite-gap integration is applicable to Heun's differential equation for the case $\gamma , \delta , \epsilon , \alpha - \beta \in \Zint +1/2$ (\cite{TV,GW,Smi,TakR}), and monodromy of Heun's differential equation can be calculated in terms of hyperelliptic integral (\cite{Tak3}) and by Hermite-Kricherver Ansats (\cite{Tak4}).

We now change variables of Heun's differential equation.
Let $\wp (x)$ be the Weierstrass doubly-elliptic function, which has double-periodicity $\wp(x)=\wp (x+2\omega _1)=\wp (x+2\omega _3)$.
We set $\omega _0=0$, $\omega _2=-\omega_1 -\omega _3$, $e_i= \wp (\omega _i)$ $(i=1,2,3)$ and 
\begin{eqnarray}
& z=\frac{\wp (x) -e_1}{e_2-e_1}, \; t=\frac{e_3-e_1}{e_2-e_1}, \;  y\tilde{\Phi}(z)=f(x), \quad \tilde{\Phi}(z)=z^{\frac{-l_0}{2}}(z-1)^{\frac{-l_1}{2}}(z-t)^{\frac{-l_2}{2}}. 
\end{eqnarray}
Then Heun's differential equation (\ref{eq:Heun}) is transformed to 
\begin{equation}
H^{(l_0,l_1,l_2,l_3)} f(x)= E f(x) ,\quad H^{(l_0,l_1,l_2,l_3)}= -\frac{d^2}{dx^2} + \sum_{i=0}^3 l_i(l_i+1)\wp (x+\omega_i),
\label{eq:ellHeun}
\end{equation}
where 
\begin{eqnarray}
& l_0= \beta -\alpha -1/2,\; l_1= -\gamma +1/2, \; l_2=-\delta +1/2, \; l_3=-\epsilon +1/2. 
\end{eqnarray}
The singularities $\{ 0 ,1,t,\infty \}$ of Eq.(\ref{eq:Heun}) correspond to half-periods $\{ \omega _1, \omega _2, \omega _3, \omega _0  \}$ of the elliptic 
function $\wp (x)$.
The parameter $q$ essentially corresponds to the eigenvalue $E$ (see \cite{Tak2}). 
Note that the operator $H^{(l_0,l_1,l_2,l_3)} $ is Hamiltonian of Inozemtsev system of type $BC_1$.
Lam\'e's  differential equation is obtained as a special case $l_1=l_2=l_3=0$ ($\gamma =\delta = \epsilon =1/2$).

By the way, Darboux transformation is extensively applied to soliton theory (\cite{MS}), and one of the features is isospectrality.
Here we briefly review Darboux transformation for Schr\"odinger equations, which is simpler than that for soliton equations.

Set 
$H = -\frac{d^2}{dx^2} + q(x)$ and assume that $\phi _0(x)$ is an eigenfunction, i.e. $H\phi _0(x) =E_0 \phi _0(x)$ for some $E_0$.
Then $q(x)=\phi _0''(x)/\phi _0(x)+E_0$.
Define an annihilation operator and a creation operator by $L=\frac{d}{dx} -\frac{\phi _0'(x)}{\phi _0(x)}$, $L^{\dagger}=-\frac{d}{dx} -\frac{\phi _0'(x)}{\phi _0(x)}$.
Then $L^{\dagger }L  = -\frac{d^2}{dx^2} +\left( \frac{\phi _0'(x)}{\phi _0(x)}\right) ' + \left( \frac{\phi _0'(x)}{\phi _0(x)}\right) ^2 = H-E_0$.
By exchanging $L^{\dagger }$ and $L $, we have $ LL^{\dagger }= -\frac{d^2}{dx^2} -\left( \frac{\phi _0'(x)}{\phi _0(x)}\right) ' + \left( \frac{\phi _0'(x)}{\phi _0(x)}\right) ^2 = H-2\left( \frac{\phi _0'(x)}{\phi _0(x)}\right) ' -E_0$.
Hence we obtain the Schr\"odinger operator $\tilde{H}= -\frac{d^2}{dx^2} +q(x)+ 2\left( \frac{\phi _0'(x)}{\phi _0(x)}\right) '$ and the relations $H=L^{\dagger }L+E_0$, $\tilde{H}=LL^{\dagger }+E_0$.
It follows from the relation $ LH=LL^{\dagger }L+E_0L=\tilde{H}L$ that, if $f(x)$ is an eigenfunction of $H$ with the eigenvalue $E$, then $Lf(x)$ is an eigenfunction of $\tilde{H}$ with the eigenvalue $E$, because $\tilde{H}(Lf(x))=LHf(x)=L(Ef(x))=E(Lf(x))$.
The transformation from $H$ to $\tilde{H}$ (or the operation by $L$) ia called Darboux transformation.
Note that the operator $L$ annihilates the 1-dimensional space ${\bf C}\phi_0(x)$.

Darboux transformation is applied to isomonodromic property of Heun's differential equations, which we will explain by an example.
Let $H^{(2,0,0,0)}$ be the operator defined in Eq.(\ref{eq:ellHeun}).
It follows from a straightforward calculation that the function $\phi _0 (x) = \sqrt{(\wp (x) -e_1) (\wp (x)-e_2)}$ is an eigenfunction of $H^{(2,0,0,0)} $ with the eigenvalue $3e_3$, i.e. $ H^{(2,0,0,0)} \phi _0 (x) = 3 e_3 \phi _0 (x)$.
The annihilation operator and the creation operator are expressed as $L=\frac{d}{dx} -\frac{\wp ' (x)}{2(\wp(x)-e_1)} -\frac{\wp ' (x)}{2(\wp (x)-e_2 )}$, $L^{\dagger}=-\frac{d}{dx} -\frac{\wp ' (x)}{2(\wp(x)-e_1)} -\frac{\wp ' (x)}{2(\wp (x)-e_2 )}$ and we have $\tilde{H}= H + 2\left( \frac{\phi _0'(x)}{\phi _0(x)}\right) ' = H^{(1,1,1,0)} $, $ H^{(1,1,1,0)}  L =L  H^{(2,0,0,0)}$.
The operator $L$ determines a map from an eigenfunction of $H^{(2,0,0,0)} $ with an eigenvalue $E$ to the eigenfunction of $H^{(2,0,0,0)} $ with the eigenvalue $E$ for all $E$.
Since the operator $L$ is doubly-periodic, the monodromy of solutions to $(H^{(2,0,0,0)} -E)f(z)=0$ coincides with the one to $(H^{(1,1,1,0)} -E)f(z)=0$.
In other words, the operator $L = \frac{d}{dx} -\frac{\wp ' (x)}{2(\wp(x)-e_1)} -\frac{\wp ' (x)}{2(\wp (x)-e_2 )}$ induces isomonodromic property of $H^{(2,0,0,0)} $ and $H^{(1,1,1,0)}$.
Note that monodromy of solutions to $(H^{(2,0,0,0)} -E)f(z)=0$ (resp. $(H^{(1,1,1,0)} -E)f(z)=0$) had been calculated individually by finite-gap integration (\cite{Tak4}), and it was found later that Darboux transformation provides a direct connection of coincidence of monodromy (\cite{KS0,Tak5}).
Calculation of monodromy is important for application to physics, because monodromy is used judge periodicity of eigenfunctions and periodicity is related to boundary condition for the eigenfunctions.
In this paper we generalize Darboux transformation and apply it to isomonodromic property for Heun's differential equation by following \cite{Tak5}.

By the way, an integral transformation on Heun's differential equation was noticed in a connection with middle convolution (\cite{TakI,TakM}), although it had been established before by Kazakov and Slavyanov \cite{KS}.
In this paper we will study the integral transformation on Heun's differential equation and observe relationship with generalized Darboux transformation.
\section{Darboux-Crum transformation}
We introduce a proposition that is a generalization of Darboux transformation.
\begin{prop} (\cite{AST,Tak5,Crum}) \label{prop:GDTA}
Let $H= -\frac{d^2}{dx^2} + q(x)$ and $U$ be a $n$-dimensional space of functions which is invariant under the action of $H$.
Let 
\begin{eqnarray}
& L= \left( \frac{d}{dx} \right) ^n +\sum _{i=1}^n c_i(x) \left( \frac{d}{dx} \right) ^{n-i}
\label{eq:L}
\end{eqnarray}
be the operator that annihilates any elements in $U$, i.e. $Lf(x)=0$ for all $f(x) \in U$.
Set $\tilde{H}= -\frac{d^2}{dx^2} +q(x)+ 2c'_1(x)$.
 then we have
\begin{equation}
\tilde{H}L=LH.
\end{equation}
\end{prop}
We call $L$ the generalized Darboux transformation or Darboux-Crum transformation.
Note that Crum \cite{Crum} had obtained a proposition whose expression is slightly different.
If $n=1$, then we reproduce Darboux transformation.

In order to apply Proposition \ref{prop:GDTA} for Heun's differential equation, we recall quasi-solvability of Heun's equation.
\begin{prop} \label{findim} (\cite[Proposition 5.1]{Tak2})
Let $\alpha _i$ be a number such that $\alpha _i= -l_i$ or $\alpha _i= l_i+1$ for each $i\in \{ 0,1,2,3\} $.
Set $d=-\sum_{i=0}^3 \alpha _i /2$ and assume $d\in \Zint_{\geq 0}$.
Let $V_{\alpha _0, \alpha _1, \alpha _2, \alpha _3}$ be the $d+1$-dimensional space spanned by 
$ \left\{ \widehat{\Phi}(\wp (x)) \wp(x)^n\right\} _{n=0, \dots ,d},$
where $\widehat{\Phi}(z)=(z-e_1)^{\alpha _1/2}(z-e_2)^{\alpha _2/2}(z-e_3)^{\alpha _3/2}$.
Then the operator $H^{(l_0,l_1,l_2, l_3)}$ preserves the space $V_{\alpha _0, \alpha _1, \alpha _2, \alpha _3}$.
\end{prop}
Let $ L_{\alpha _0, \alpha _1, \alpha _2, \alpha _3} $ be the monic differential operator of order $d+1$ which annihilates the space $V_{\alpha _0, \alpha _1, \alpha _2, \alpha _3}$.
By calculating the coefficient $c_1 (x)$ in Eq.(\ref{eq:L}) for the operator $ L_{\alpha _0, \alpha _1, \alpha _2, \alpha _3} $, we have Darboux-Crum transformation for Heun's differential equation.
\begin{thm} \label{thm:HL0123L0123H} (\cite[Theorem 3.3]{Tak5})
Let $\alpha_i$ be a number such that $\alpha_i= -l_i$ or $\alpha_i= l_i+1$ for each $i\in \{ 0,1,2,3\} $. Set $d=-\sum_{i=0}^3 \alpha_i /2$ and assume $d\in \Zint_{\geq 0}$. Then we have
\begin{eqnarray}
& H^{(\alpha _0 +d,\alpha _1 +d,\alpha _2 +d,\alpha _3 +d)}  L_{\alpha _0, \alpha _1, \alpha _2, \alpha _3}  =L _{\alpha _0, \alpha _1, \alpha _2, \alpha _3} H^{(l_0,l_1,l_2,l_3)}.
\end{eqnarray}
\end{thm}

Let $f_1(x,E)$, $f_2 (x,E)$ be a basis of solutions to $(H^{(l_0,l_1,l_2,l_3)} -E) f(x)=0$.
Since the operator $H^{(l_0,l_1,l_2,l_3)}$ is doubly-periodic, the functions 
$f_1(x+2\omega _k,E)$, $f_2(x+2\omega _k,E)$ $(k=1,3)$ are also solutions to the differential equation.
Let $M_{2\omega _k} (E)$ $(k=1,3)$ be a monodromy matrix on the shift $x \rightarrow x+ 2\omega _k$ with respect to the basis $\{ f_ 1(x, E)$, $f_2 (x ,E) \}$, i.e.
\begin{equation}
( f_1(x+2\omega _k,E) \; f_2(x+2\omega _k,E)) = (f_1(x,E) \; f_2(x,E))  M_{2\omega _k} (E) .
\end{equation}
It follows from absence of first order derivative of the differential equation that $\det M_{2\omega _k} (E) =1$ $(k=1,3)$.
Note that ${\rm{tr} } M_{2\omega _k} (E) $ is independent from the choices of a basis of solutions to $(H^{(l_0,l_1,l_2,l_3)} -E) f(x)=0$.
Set $\tilde{f}_i(x,E) =L_{\alpha _0, \alpha _1, \alpha _2, \alpha _3}f_i(x,E)$ $(i=1,2)$. 
Then
$H^{(\alpha _0 +d,\alpha _1 +d,\alpha _2 +d,\alpha _3 +d)}\tilde{f}_i(x,E) =E\tilde{f}_i(x,E)$.
Since $L_{\alpha _0, \alpha _1, \alpha _2, \alpha _3}$ is doubly-periodic, we have
\begin{eqnarray}
& (\tilde{f}_1(x+2\omega _k,E) \; \tilde{f}_2(x+2\omega _k,E))=(\tilde{f}_1(x,E) \; \tilde{f}_2(x,E)) M_{2\omega _k} (E) .
\end{eqnarray}
Hence the monodromy structure of $H^{(l_0,l_1,l_2,l_3)}$ with respect to a shift of a period coincides with the one of $H^{(\alpha _0 +d,\alpha _1 +d,\alpha _2 +d,\alpha _3 +d)}$,
and the operator $L_{\alpha _0, \alpha _1, \alpha _2,\alpha _3}$ defines an isomonodromic transformation from $H^{(l_0,l_1,l_2,l_3)}$ to $H^{(\alpha _0 +d,\alpha _1 +d,\alpha _2 +d,\alpha _3 +d)}$.
In particular, we have ${\rm{tr}} M_{2\omega _k} (E) = {\rm{tr}} \tilde{M}_{2\omega _k} (E) $, where ${\rm{tr}} M_{2\omega _k} (E)$ (resp. ${\rm{tr}} \tilde{M}_{2\omega _k} (E)$) is the trace of the monodromy matrix of solutions to $(H^{(l_0,l_1,l_2,l_3)} -E) f(x)=0$ (resp. $(H^{(\alpha _0 +d,\alpha _1 +d,\alpha _2 +d,\alpha _3 +d)} -E) f(x)=0$).

We now provide an example of Darboux-Crum transformation.
On the case $l_0=2l$ $(l\in \Zint_{\geq 1})$, $l_1=l_2=l_3=0$, we set $\alpha _0=-2l$, $\alpha _1=\alpha _2=1$, $\alpha _3=0$.
Then $d=-(\alpha _0+\dots +\alpha _3)/2=l-1$ and we may adopt Theorem \ref{thm:HL0123L0123H}.
We have
\begin{eqnarray}
& H^{(l,l,l,l-1)}  L_{-2l, 1, 1, 0}  =L _{-2l,1,1,0} H^{(2l,0,0,0)},
\end{eqnarray}
because $ H^{(-l-1,l,l,l-1)}=H^{(l,l,l,l-1)}$.
Hence $H^{(2l,0,0,0)}$ is isomonodromic to $H^{(l,l,l,l-1)}$ for $l \in \Zint _{\geq 1}$.
If $l=1$, then $d=0$, $ H^{(1,1,1,0)}  L_{-2, 1, 1, 0} =L _{-2,1,1,0} H^{(2,0,0,0)}$,
and the operator $L_{-2, 1, 1, 0}$ is written as 
$L_{-2,1,1,0}=\frac{d}{dx} -\frac{\wp ' (x)}{2(\wp(x)-e_1)} -\frac{\wp ' (x)}{2(\wp (x)-e_2 )}$.
Namely we reproduce the example in introduction.

\section{Application to finite-gap integration}

The method of finite-gap integration has been studied more than 30 years (\cite{Mat}) and a target of earlier study was KdV equation.
It was also applied to Heun's differential equation (\cite{TV,GW,Smi,TakR}) with the condition $l_0, l_1, l_2, l_3 \in \Zint  $ and to a certain class of Fuchsian differential equations (\cite{TakZ}).
Consequently we have formulas of monodnomy of Heun's differential equation by hyperelliptic integrals (\cite{Tak3}), Bethe Ansats (\cite{Tak1}) and Hermite-Krichever Ansatz (\cite{Tak4}).
A remarkable feature of finite-gap integration is existence of an odd-order differential operator $\tilde{A}$ which commutes with the operator $H^{(l_0,l_1,l_2,l_3)}$, i.e. $[\tilde{A}, H^{(l_0,l_1,l_2,l_3)} ] =0$.

We can construct the odd-order differential operator $\tilde{A}$ by composing Darboux-Crum transformations for the case $l_0, l_1, l_2, l_3 \in \Zint $.
\begin{prop} (\cite{Tak5})
If $l_0, l_1, l_2, l_3 \in \Zint $, then we can construct an odd-order differential operator $\tilde{A}$
such that $[\tilde{A}, H^{(l_0,l_1,l_2,l_3)} ] =0$ by composing four Darboux-Crum transformations.
\end{prop}
For details, see \cite[Proposition 6.1]{Tak5}.
If $l_0=2$, $l_1=l_2=l_3=0$, then there exists a differential operator $\tilde{A}$ of fifth order written as 
$\tilde{A} = L_{2,-1,-1,0} L_{1,-2,1,0} L_{0,2,-1,-1}  L_{-2,0,0,0} $ which satisfies $[\tilde{A}, H^{(2,0,0,0)} ] =0$.

\section{Integral transformation of Heun's differential equation}
We introduce integral representation of Heun's differential equation as a generalization of Darboux transformation.

Integral transformation on Heun's differential equation has been noticed in a connection with middle convolution (\cite{TakI,TakM}), although it had been established before by Kazakov and Slavyanov \cite{KS}.
In order to formulate the integral transformation, we fix a base point $o $ of the integrals in the complex plane $\Cplx$ appropriately.
Let $p $ be an element of the Riemann sphere $\Cplx \cup \{\infty \}$ and $\gamma _p$ be a cycle in the Riemann sphere with the variable $w$ which starts from $w=o$, turns the point $w=p$ anti-clockwise and ends at $w=o$.
Let $[\gamma _z ,\gamma _p] = \gamma _z \gamma _p \gamma _z ^{-1} \gamma _p ^{-1}$ be the Pochhammer contour.
\begin{prop} (\cite{KS,TakM}) \label{prop:Heunint}
(i) Set
\begin{eqnarray}
& \{\mu -(2-\alpha )\}\{\mu -(2-\beta )\}=0 , \; \gamma '=\gamma +\mu -1, \; \delta' =\delta +\mu -1, \; \epsilon '=\epsilon +\mu -1,\nonumber \\
& \alpha '=\mu , \; \beta '= 2\mu + \alpha +\beta -3 , \; q'=q+(1-\mu )(\epsilon +\delta t+(\gamma -\mu ) (t+1)). 
\end{eqnarray}
Let $y(w)$ be a solution to 
\begin{eqnarray}
& \frac{d^2y}{dw^2} + \left( \frac{\gamma }{w}+\frac{\delta }{w-1}+\frac{\epsilon }{w-t}\right) \frac{dy}{dw} +\frac{\alpha  \beta  w -q}{w(w-1)(w-t)} y=0.
\end{eqnarray}
Then the functions $(i \in \{ 0,1,t,\infty \})$ 
\begin{eqnarray}
& \tilde{y}(z)=\int _{[\gamma _z ,\gamma _i]} y(w) (z-w)^{-\mu } dw 
\end{eqnarray}
are solutions to 
\begin{eqnarray}
& \frac{d^2\tilde{y}}{dz^2} + \left( \frac{\gamma '}{z}+\frac{\delta '}{z-1}+\frac{\epsilon '}{z-t}\right) \frac{d\tilde{y}}{dz} +\frac{\alpha '\beta 'z -q'}{z(z-1)(z-t)} \tilde{y}=0.
\label{eq:Heunprime}
\end{eqnarray}
(ii) Under the notation of (i), we assume $\mu \in \Zint _{\geq 1}$ additionally.
Then
\begin{eqnarray}
& \tilde{y}(z)= \left( \frac{d}{dz} \right) ^{\mu -1 }y(z)
\end{eqnarray}
is a solution to Eq.(\ref{eq:Heunprime}). 
\end{prop}

We rewrite Proposition \ref{prop:Heunint} to the form of the elliptical representation by setting $\mu =d+2$.
It is remarkable that the eigenvalue $E$ is unchanged by the integral transformation.
\begin{prop} (\cite{TakIT}) \label{prop:ellipinttras}
(i) Let $\sigma (x)$ be the Weierstrass sigma function, $\sigma _i (x)$ $(i=1,2,3)$ be the Weierstrass co-sigma function which has a zero at $x=\omega _i$, and $I_i$ $(i=0,1,2,3)$ be the cycle on the complex plane with the variable $y$ such that points $y=x$ and $y=-x+2\omega _i$ are contained and the half-periods $\Zint \omega _1 +\Zint \omega _3$ are not contained inside the cycle.
Let $\alpha _i$ be a number such that $\alpha _i= -l_i$ or $\alpha _i= l_i+1$ for each $i\in \{ 0,1,2,3\} $.
Set $d=-\sum_{i=0}^3 \alpha _i /2$.
Let $f(x)$ be a solution to $(H^{(l_0,l_1,l_2,l_3)} -E) f(x)=0$.
Then the functions 
\begin{eqnarray}
& \tilde{f}(x)=\sigma (x) ^{\alpha _0 +d+1 } \sigma _1 (x) ^{\alpha _1 +d+1 } \sigma _2 (x) ^{\alpha _2 +d+1} \sigma _3 (x) ^{\alpha _3 +d+1} \cdot \label{eq:Hintell2} \\
& \quad \quad \int _{I_i } f (y) \sigma (y) ^{1 -\alpha _0}  \sigma _1 (y) ^{1 -\alpha _1} \sigma _2 (y) ^{1 -\alpha _2} \sigma _3 (y) ^{1 -\alpha _3} (\sigma (x+y) \sigma (x-y) )^{-d-2 } dy \nonumber 
\end{eqnarray}
$(i \in \{ 0,1,2,3 \})$ are solutions to $(H^{(\alpha _0 +d ,\alpha _1 +d ,\alpha _2 +d ,\alpha _3 +d )} -E) f(x)=0$. 

(ii) Under the notation of (i), we assume $d\in \Zint_{\geq -1}$ additionally.
Let $f(x)$ be a solution to $(H^{(l_0,l_1,l_2,l_3)} -E) f(x)=0$.
Then the function
\begin{equation}
\tilde{f}(x)=\wp '(x)^{d+1} \prod _{i=1}^3 (\wp (x) -e_i)^{\alpha _i /2} 
\left( \frac{1}{\wp '(x) } \frac{d}{dx} \right) ^{d+1} \left\{ f(x) \prod _{i=1}^3 (\wp (x) -e_i)^{-\alpha _i /2} \right\} 
\label{eq:tfdiff}
\end{equation}
is a solution to $(H^{(\alpha _0 +d ,\alpha _1 +d ,\alpha _2 +d ,\alpha _3 +d )} -E) f(x)=0$. 
\end{prop}
If $d\in \Zint_{\geq 0}$ and $f(x) \in V_{\alpha _0, \alpha _1, \alpha _2, \alpha _3}$, then the function $\tilde{f}(x)$ in Eq.(\ref{eq:tfdiff}) is identically equal to zero, and it follows that 
\begin{equation}
\wp '(x)^{d+1} \prod _{i=1}^3 (\wp (x) -e_i)^{\alpha _i /2} \circ
\left( \frac{1}{\wp '(x) } \frac{d}{dx} \right) ^{d+1} \circ \prod _{i=1}^3 (\wp (x) -e_i)^{-\alpha _i /2} = L_{\alpha _0, \alpha _1, \alpha _2, \alpha _3}.
\label{eq:tfdiffL}
\end{equation}
Therefore we may regard the integral transformation as  a generalization of Darboux-Crum transformation by removing the condition $d\in \Zint_{\geq 0}$.

We can express the monodromy of solutions to $(H^{(\alpha _0 +d ,\alpha _1 +d ,\alpha _2 +d ,\alpha _3 +d )}-E) f(x)=0$ in terms of the monodromy of solutions to $(H^{(l_0,l_1,l_2,l_3)} -E) f(x)=0$ by applying integral transformations, which was performed in \cite{TakIT}.
\begin{thm} (\cite{TakIT}) \label{thm:pp}
Let $\alpha _i$ be a number such that $\alpha _i= -l_i$ or $\alpha _i= l_i+1$ for each $i\in \{ 0,1,2,3\} $.
Set $d=-\sum_{i=0}^3 \alpha _i /2$.
Let $k \in \{1,3\}$ and $M_{2\omega _k} (E)$ (resp. $\tilde{M}_{2\omega _k} (E)$) be the monodromy matrix by the shift of the period $x \rightarrow x+2\omega _k$ with respect to a certain basis of solutions to $(H^{(l_0,l_1,l_2,l_3)} -E) f(x)=0$ (resp. $(H^{(\alpha _0 +d ,\alpha _1 +d ,\alpha _2 +d ,\alpha _3 +d )}-E) f(x)=0$).
Then $\mbox{\rm{tr}}M_{2\omega _k} (E)= \mbox{\rm{tr}}\tilde{M} _{2\omega _k}(E)$.
\end{thm}
Hence the monodromy structure of $H^{(l_0,l_1,l_2,l_3)}$ with respect to a shift of a period coincides with the one of $H^{(\alpha _0 +d,\alpha _1 +d,\alpha _2 +d,\alpha _3 +d)}$ for all $d(= -\sum_{i=0}^3 \alpha _i /2 )$.
As a corollary we have 
\begin{cor} (\cite{TakIT}) \label{cor:pp}
We keep the notations in Theorem \ref{thm:pp}.
Let $k \in \{1,3\}$.
If there exists a non-zero solution $f(x,E)$ to $(H^{(l_0,l_1,l_2,l_3)} -E) f(x,E)=0$ such that $f(x+2\omega _k,E) =C_k(E) f(x,E) $, then there exists a non-zero solution $\tilde{f}(x,E)$ to $(H^{(\alpha _0 +d ,\alpha _1 +d ,\alpha _2 +d ,\alpha _3 +d )}-E) \tilde{f}(x,E)=0$ such that $\tilde{f}(x+2\omega _k,E) =C_k(E) \tilde{f}(x,E)$.
In other word, periodicity is preserved by the integral transformation.
\end{cor}
For example the monodromy structure of $H^{(2l,0,0,0)}$ with respect to a shift of a period coincides with the one of $H^{(l,l,l,l-1)}$ for all $l\in \Rea$.
In particular, the condition of eigenvalues $E$ such that there exists a non-zero periodic solution to $H^{(2l,0,0,0)} f(x) = Ef(x) $ with respect to the shift $x \rightarrow x+2\omega _1$ coincides with that of eigenvalues $E$ such that there exists a non-zero periodic solution to $H^{(l,l,l,l-1)} f(x) = Ef(x)$ with respect to the shift $x \rightarrow x+2\omega _1$ for all $l\in \Rea$.

\section{Concluding remarks}
In this paper we observed some aspects of Darboux-Crum transformation and integral transformation on Heun's differential equation.
In particular we can regard Darboux-Crum transformation as a specialization of integral transformation, and periodicities of solutions to Heun's differential equation is preserved by the transformations.
We may interpret our result as a different realization of the observation by Khare and Sukhatme \cite{KS0}.

We note that integral transformation may produce new solutions to Heun's differential equation.
For example, solutions of Heun's differential equation with the condition $\gamma ,\delta , \epsilon , \alpha +1/2, \beta +1/2 \in \Zint $ (i.e. $l_0, l_1, l_2, l_3 \in \Zint +1/2$, $l_0+l_1+l_2+l_3 \in 2\Zint +1$) can be expressed in use of finite-gap solutions,
and the special case $\gamma = \delta = \epsilon =1$, $\alpha =1/2$, $\beta =3/2 $ (i.e. $l_0=1/2$, $l_1=l_2=l_3=-1/2$) was previously studied by Valent \cite{Val} to understand an eigenvalue problem related to certain birth and death processes.
We hope further applications of new solutions to physics and mathematics.

\end{document}